\documentclass[preprint,12pt]{elsarticle}

	\usepackage{graphicx}
\usepackage{amssymb}
\usepackage{amsthm}
\usepackage{url}
\usepackage{lineno}
\usepackage{amssymb}
\usepackage{amsmath}
\usepackage{makeidx}
\usepackage{amscd}
\usepackage{fancyhdr}
\usepackage{yfonts}
\usepackage{graphicx}
\usepackage{bbm}
\usepackage{enumitem}
\usepackage{amsthm}
\usepackage{listings}
\usepackage{fancyvrb}
\usepackage{url}
\usepackage{mathrsfs}  
\usepackage{color}

\journal{Finite fields and their applications}

\newtheorem{theo.}{Theorem}[section]
\newtheorem{lem.}[theo.]{Lemma}
\newtheorem{cor.}[theo.]{Corollary}
\newtheorem{prop.}[theo.]{Proposition}
\newtheorem{def.}[theo.]{Definition}
\newtheorem{ex.}[theo.]{Example}
\newtheorem{not.}[theo.]{Notation}

\newcommand{\F}{\ensuremath{\mathbb F}}
\newcommand{\Z}{\ensuremath{\mathbb Z}}
\newcommand{\HC}{\ensuremath{\mathcal H}}

\newcommand{\contradiction}{{\hbox{%
    \setbox0=\hbox{$\mkern-3mu\times\mkern-3mu$}%
    \setbox1=\hbox to0pt{\hss$\times$\hss}%
    \copy0\raisebox{0.5\wd0}{\copy1}\raisebox{-0.5\wd0}{\box1}\box0
}}}

\begin{document}

\begin{frontmatter}


\title{Geometric primality tests using curves of genus 1 \& 2}



\author{Eduardo Ru\'{i}z Duarte}

\address{University of Groningen, The Netherlands}

\begin{abstract}
We revisit and generalize some geometric techniques behind deterministic primality testing for some integer sequences using curves of genus $1$ over finite rings. Subsequently we develop a similar primality test using the Jacobian of a genus $2$ curve. 
\end{abstract}

\begin{keyword}
Primality \sep Elliptic Curves \sep Hyperelliptic Curves


\end{keyword}

\end{frontmatter}


\section{Introduction}
\label{S:1}
\newcommand{\ja}{\mathcal{J}}

\newcommand{\hc}{\mathcal{H}}

This paper is mainly inspired by a lecture at the \textit{Intercity Seminar} \cite{intercitytop} given by Jaap Top, \textit{``Lucas-Lehmer revisited"}. Also a short note by Dick Gross was relevant for this topic since he produced the first deterministic primality test for numbers of the form $2^p-1$ (Mersenne numbers) using an elliptic curve \cite{grossmersenne}. We begin with the simplest case to explore and generalize the usage of the elliptic curve $E_t:y^2=x^3 - (t^2 + 1)x$ for primality testing of integers $m2^n-1$ using the $\Z$-module structure of $E_t$.   Further, in \cite{densav} Denomme and Savin used complex multiplication to develop several primality tests for different sequences of integers, later some of them generalized by Gurevich and Kunyavski{\u\i} in \cite{gurkun}. Particularly, Denomme and Savin used $E:y^2 = x^3 - x$ to do a primality test on Fermat numbers using the $\text{End}_{\mathbb{Q}(i)}(E)$-module structure of $E$. Here we revisit and extend their setting from Fermat integers to integers of the form $p^2 16^n + 1$ where $p\equiv \pm 1\bmod 10$ and $p<2^n$. Furthermore, with this, we answer an open question stated in \cite{primeframework} by Abatzoglou, Silverberg, Sutherland and Wong (see Remark 4.13). 
This question asks about the design of a potential primality test using Jacobians of genus $2$ curves. This new primality test is designed for integers of the form $4\cdot 5^{n}-1$ and it uses the Jacobian $\ja$ of the genus $2$ curve $\hc:y^2=x^5 + h$ as a cyclic $\text{End}_{\mathbb{Q}(\sqrt{5})}(\ja)$-module. We emphasize that more efficient primality tests for integers $4\cdot 5^n -1$ may exist, but here we state a theoretical result to do a primality test for these integers using an Abelian variety of dimension 2.      \\ 

\subsection{Primality testing \textit{\`a la} Lucas}
\noindent
It is well known that a necessary condition (but not sufficient) for a number $n\in\mathbb{N}$ 
to be prime is that for all $a\in\mathbb{N}$ such that $2\leq a<n$ the congruence $a^{n-1}\equiv 
1\bmod n$ holds. This \textit{Little Theorem} by Fermat can be used as a \textit{test} for compositeness calculating the congruence for several $a$. We infer that $n$ is composite if for some $a$, the congruence does not hold. When the congruence holds for many choices of $a$ the number $n$ is said to be \textit{probably prime}. The computation of this congruence can be done quite fast using modular repeated squaring.\\ \\
Unfortunately there is a problem with this \textit{Fermat test}, there are infinitely many composite numbers such as $m=561$ satisfying $a^{m-1}\equiv 1\bmod m$ for all $a$ such that $(m,a)=1$. These numbers are known as \textit{Carmichael numbers}. Even though Carmichael numbers are \text{rarer} than prime numbers (see \cite{erdcar}) other extensions of this test were developed to deal with this, like Miller-Rabin or Solovay-Strassen which are more common in practice. In order to turn this \textit{Fermat test} into a primality testing algorithm, \'{E}douard Lucas stated the following theorem:
\begin{theo.}\label{lucast}[Lucas, 1876] Let $a,n\in\mathbb{Z}$ such that $a^{n-1}\equiv 1\bmod n$ and $a^{\frac{n-1}{p}}\not\equiv 1\bmod n$ for all primes $p\mid(n-1)$. Then $n$ is prime.
\end{theo.}
\begin{proof}
Let $a\in{(\mathbb{Z}/n\mathbb{Z}})^{\times}$ and $k=\#\langle a\rangle$. Since $a^{n-1}\equiv 1\bmod n$ then $k\mid n-1$. Further, we have that $a^{\frac{n-1}{p}}\not \equiv 1 \bmod n$ for all $p\mid (n-1)$, hence $k=n-1$. With this we have that $\#{(\mathbb{Z}/n\mathbb{Z}})^{\times}= n-1$ and then $n$ is prime. 
\end{proof}
\noindent
This elementary theorem is used by several deterministic primality tests.  \\
A problem for potential algorithms that could arise from this theorem is that it requires the prime divisors of $n-1$. This is very difficult in general, but for example, if we restrict our algorithms to potential prime numbers of the form $k2^n+1$ or $2^{2^n}+1$ this theorem can be applied effectively. Also, another less difficult problem when using Theorem~\ref{lucast} is that in case of $n$ being prime, we need to find a correct $a\in(\mathbb{Z}/n\mathbb{Z})^{\times}$ that satisfies the hypotheses of Theorem \ref{lucast}. This ``\textit{problem}'' means that when $n$ is prime then $(\mathbb{Z}/n\mathbb{Z})^{\times}$ is cyclic, so we need an $a\in (\mathbb{Z}/n\mathbb{Z})^{\times}$ that generates this cyclic group (of units of $\mathbb{Z}/n\mathbb{Z}$). The existence of this $a$ satisfying Theorem \ref{lucast} is a classical result by Gau\ss\space exposed in the article 57 from \textit{Disquisitiones} where he calls them \textit{primitive roots modulo $n$}, see  \cite{gaussdisq}. Gau\ss\space proved that such $a$ exists but finding it in general is a different problem.\\
The main motivation for primality testing using geometrical tools can be traduced from the following theorem for Mersenne numbers $M_n:=2^n-1$.
\begin{theo.}\label{lucaslehmer}[Lucas-Lehmer]
Consider the sequence $a_0:=4, a_{i+1}:=a_i^2 -2$. Let $n>2$. $M_n:=2^n-1$ is prime if and only if $a_{n-2}\equiv 0\bmod M_n$.
\end{theo.}
We sketch the proof of this theorem using the properties of the Abelian group given by a Pell conic, namely: 
\[
G(k):=\{ (x,y)\in k\times k: x^2-3y^2=1\}.
\]
The element $a_j$ in the sequence of Theorem \ref{lucaslehmer} is exactly $2$ times the $x$ coordinate of the $j^{th}$ recursive squaring of the point $(2,1)\in G\bmod M_n$, namely $a_j=2\cdot x(2^j(2,1)) \bmod M_n$. When $M_n$ is prime, it can be proved that the point $(2,1)$ is not twice some other point in $G(\F_{M_n})\cong\Z/(2^n)$ and $a_{n-2}$  corresponds to a point of order $4$ in $G$ given by $2^{n-2}(2,1)=(0,\pm\tfrac{1}{\sqrt{3}})$. \\ When the sequence holds and one assumes that $M_n$ is not prime, a contradiction arises by taking a non-trivial prime divisor $k|M_n$. More precisely, you will encounter the inequality $2^n\leq \#G(\F_k)=k \pm 1$, where $2^n$ is the size of the subgroup of $G(\F_k)$ generated by the point $(2,1)$. 

In the following sections we will explore a geometrical perspective of the \textit{Lucasian} primality tests arising from these ideas using algebraic groups from  elliptic curves and from Jacobians of genus $2$ curves.
\newcommand{\mcyc}{\mathbb{G}_{m,\mathbb{Z}}}
\newcommand{\fer}{\mathcal{F}}
\newcommand{\pr}{\mathcal{A}}
\newcommand{\legendre}[2]{\genfrac{(}{)}{}{}{#1}{#2}}
\section{Primality testing with genus 1 curves}
In this section we construct a primality test using properties of supersingular elliptic curves without using complex multiplication; later we will use
complex multiplication as well.\\
For the first part we will use recursive doubling of points similar to the primality test algorithm proposed by Dick Gross for Mersenne primes but now for integers of the form $m2^n -1$. \\
Additionally we will extend a test presented by Denomme and Savin in \cite{densav} from Fermat numbers to integers of the form $p^2 16^n + 1$ where $p\equiv\pm 1\bmod 10$ and $p<4^n$. The idea behind the test by Denomme and Savin is to use an endomorphism of degree $2$ arising from the complex multiplication of an elliptic curve $E$ of $j$-invariant $1728$, namely  $(1+i)\in \text{End}(E)$. Their method is to recursively apply this map on a specific point to prove that a Fermat number is prime using the same principle given by Theorem \ref{lucaslehmer}. They use the $\Z[i]$-module structure of the elliptic curve $E$ gotten by the action of $\Z[i]$ on the Abelian group given by the rational points of the elliptic curve $E$. 
\subsection{Primality testing with supersingular elliptic curves}
In this section we provide a family of elliptic curves that will lead to primality tests for numbers of the form $\pr_{m,n}:=m2^n-1$. The following Proposition is a key part for the design of a primality test algorithm of $\pr_{m,n}$.  

\begin{prop.}\label{elldiv} Let $p\equiv 3\bmod 4$ be a prime number and $t\in\F_p$.
The equation $y^2=x^3 -(t^2 + 1)x=f_t(x)$ over $\F_{p}$ 
defines a supersingular elliptic curve $E_t/\F_p$ and the point $(-1,t)$ is not divisible by $2$ in $E_t(\F_p)$.
\end{prop.}
\begin{proof}
The first claim can be proved directly. First, since $p\equiv 3\bmod 4$, we have that $t^2\neq -1$ for all $t\in\mathbb{F}_p$. Hence $E_t$ indeed defines an elliptic curve.\\
The fact that it is supersingular is well known, compare
\cite[V~Example~4.5]{AEC}.
For convenience we provide an alternative argument.
Let $w\in\{1,3\}$ be the number of $\F_p$-rational zeros of $f_t(x)$ and let $x_0\in\F_p$ such that  $f_t(x_0)\neq 0$. We have that $f_t(-x_0)=-f_t(x_0)$, hence using $p\equiv 3\bmod 4$, one concludes $f_t(x_0)$ is a square over $\F_p$ if and only if $f_t(-x_0)$ is not a square over $\F_p$ (this is because $-1\notin \F_p^2$). Hence, the number of points of $E_t(\F_p)$ is given by twice the aforementioned squares $f_t(x_i)$ for all $x_i\in\F_p$ such that $f_t(x_i)\neq 0$. The value $\#E_t(\F_p)$ is given by counting these $x_i$ which are $\tfrac{p-w}{2}\cdot 2$ and adding  the number of Weierstrass points given by $w+1$. Hence $\#E_t(\F_p)=p+1$ for all $t\in\mathbb{F}_p$ and $E_t/\F_p$ is supersingular. \\ \\
To prove that $(-1,t)$ is not divisible by two, in other words that there 
is no $Q\in E_t(\F_p)$ such that $2Q=(-1,t)$, consider the multiplication-
by-$2$ map given by $2\in\text{End}_{\F_p}(E)$. It is equivalent to show 
that $(-1,t)\notin 2 E_t(\F_p)$. The $2$-descent homomorphism $\delta$ (see \cite{AEC} 
Chapter X, \S 4 Prop. 4.9 for details) is 
useful here since $\text{Ker}(\delta)=2E_t(\F_p)$. We proceed to construct $\delta$ 
for $E_t(\F_p)$ and apply it to $(-1,t)$. This construction will be done in
two cases depending on $t^2+1$ being a square or not in $\F_p$.\\ \\
Let $t^2 + 1 \notin\F_p^2$ and consider the ring 
$\mathfrak{R}_{t}:=\F_p[X]/(f_t(X))$. Since $t^2+1$ is not a square, 
$f_t(X)$ defines only one affine $\F_p$-rational Weierstrass point in 
$E_t(\F_p)$, namely $(0,0)$. Hence $\mathfrak{R}_t\cong \F_p\times 
\F_p[\xi]/(\xi^2-(t^2+1))$. Let $P:=(\alpha,\beta)\in E_t(\F_p)$. Since 
$X\in\mathfrak{R}_t$ satisfies the equation $f_t(X)=0$, the $2$-descent homomorphism of 
$E_t(\F_p)$ in this case is given by:
\begin{align}\label{twodescent}
\begin{split}
\delta:E_t(\F_p)&\to\mathfrak{R}_t^\times/\mathfrak{R}_t^{\times^2}\\
P &\mapsto  
\begin{cases}
\mathfrak{1} &\text{if $P=\infty$}\\
[-(t^2+1)-X]  &\text{if $P=(0,0)$}\\
[\alpha-X]  &\text{otherwise, }\\
\end{cases}
\end{split}
\end{align}
compare \cite[\S~3.5]{ratpoints}.\\
Since $\text{Ker}(\delta)=2E_t(\F_q)$ we have that 
$(-1,t)$ is divisible by $2$ if and only if $-1-X$ is in 
$\mathfrak{R}_t^{\times^2}$.
However it is not a square since its image $-1$ 
in $\F_p^\times/{\F_p^\times}^2$ is nontrivial. 
 \\ \\
The other case is when $\lambda^2=t^2+1\in\F_p^2$, hence $f_t(X)$ splits in $\F_p[X]$ and $\mathfrak{R}_t\cong \F_p[X]/(X)\times\F_p[X]/(X+\lambda)\times\F_p[X]/(X-\lambda)\cong \F_p\times\F_p\times\F_p$. The 2-Descent map $\delta:E_t(\F_p)\to \mathfrak{R}_t^\times/{\mathfrak{R}_t^\times}^2$ in this case applied to $(-1,t)$ again yields in the first factor $-1$ which is nontrivial.
 This shows that $(-1,t)\not\in 2E_t(\F_p)$ in
 all cases. 
\end{proof}
\noindent
Now we need to know when $E_t$ is cyclic, this will depend on the base field of $E_t$. This is important in order to establish for which integers our primality testing algorithm will be useful.
\begin{lem.}\label{ellcyc}
Let $p\equiv 3\bmod 4$ be prime and $t\in\F_p$ such that $t^2+1\notin({\F_p}^{\times})^2$. Consider the elliptic curve $E_t$ given by $y^2=x^3 - (t^2 + 1)x$, then  $E_t(\F_p)$ is cyclic. 
\end{lem.}
\begin{proof}
Consider the multiplication by $p+1$ map, that is $p\!+\!1\in\text{End}(E_t)$. We have that $\text{Ker}(p+1)=E_t(\overline{\F}_p)[p\!+\!1]\cong \Z/(p\!+\!1)\Z\times\Z/(p\!+\!1)\Z$. By Proposition \ref{elldiv} $\#E_t(\F_p)=p\!+\!1$, hence $E_t(\F_p)\leq E_t(\overline{\F}_p)[p\!+\!1]$. With this, for $1\leq\alpha \leq \beta$ we have that $E_t(\F_p)\cong \Z/\alpha\Z\times\Z/\beta\Z$ such that $\alpha\mid\beta$ and $\alpha\cdot\beta=p+1$.\\
Now we show that $\alpha\mid p-1$. \\
We look at the $\alpha-$torsion. We have that $\Z/\alpha\Z\times\Z/\alpha\Z \leq \Z/\alpha\Z\times\Z/\beta\Z\cong E_t(\F_p)$, this means that $E_t(\F_p)\supset\text{Ker}(\alpha)$. Using the surjectivity of the Weil pairing (see \cite{AEC}, III, Corollary 8.1.1) there must be $P,Q\in E_t(\F_p)$ with $(P,Q)=\omega_\alpha$ where $\omega_\alpha\in\overline{\F}_p$ is a $\alpha^{th}$ root of unity. Using the fact that the Weil pairing is Galois invariant (see \cite{AEC},III, Proposition 8.1), for any $\sigma\in\text{Gal}(\F_p(\omega_\alpha)/\F_p)$ we have that $(P^\sigma,Q^\sigma)=(P,Q)^\sigma$, hence $\omega_\alpha$ is invariant under $\sigma$. This means that $\omega_\alpha\in\F_p^\times$, which implies that $\alpha\mid \#\F_p^\times$ hence $\alpha\mid p-1$.  \\
Now we have that $\alpha\mid p+1$ and $\alpha\mid p-1$ hence $\alpha\mid 2$. This means that if $\alpha=2$ then $E_t(\F_p)$ has full two-torsion. But this is not the case since $t^2 +1$ is not a square and $p\equiv 3\bmod 4$. Hence $\alpha=1$ and $E_t(\F_p)$ is cyclic.
\end{proof}
\begin{cor.}\label{ellgen}
Let $m\geq 1$ be an odd integer and suppose $\pr_{m,n}:=m2^n-1$ be prime and $n>1$. Take $t\in \F_{\pr_{m,n}}$ such that
$t^2+1$ is not a square and consider the elliptic curve $E_t/\F_{\pr_{m,n}}$. Then the point $m(-1,t)$ generates the $2$-Sylow subgroup of $E_t(\F_{\pr_{m,n}})$.
\end{cor.}
\begin{proof}
Since $n>1$ we have that $\pr_{m,n}\equiv 3\bmod 4$. Using Lemma \ref{ellcyc} and ${t^2+1}\notin {\F_{\pr_{m,n}}}^2$, the group $E_t(\F_{\pr_{m,n}})$ is cyclic and has $m2^n$ points. By Proposition \ref{elldiv} the point $(-1,t)$ is not divisible by $2$ and since $m$ is odd, $m(-1,t)$ has order $2^n$.
\end{proof}
The following simple lemma will be used to discard some trivial small divisors of $\pr_{m,n}$ for every $n>0$ and $m>2$. This lemma will make the proof of the main Theorem of this section shorter.
\begin{lem.}\label{aprcomposite}
The number $\pr_{m,n}=m2^n - 1$ is divisible by $3$ if 
and only if one of the following conditions holds:
\begin{itemize}
\item{ $m\equiv 2\bmod 3$ and $n\equiv 1\bmod 2$};
\item{ $m\equiv 1\bmod 3$ and $n\equiv 0\bmod 2$}.
\end{itemize}
Further, $\pr_{m,n}=m2^n -1$ is divisible by $5$ if 
and only if one of the following conditions holds:
\begin{itemize}
\item{ $m\equiv 1\bmod 5$ and $n\equiv 0\bmod 4$};
\item{ $m\equiv 2\bmod 5$ and $n\equiv 3\bmod 4$};
\item{ $m\equiv 3\bmod 5$ and $n\equiv 1\bmod 4$};
\item{ $m\equiv 4\bmod 5$ and $n\equiv 2\bmod 4$}.
\end{itemize}
\end{lem.}
\begin{proof}
This is clear analyzing the period of $2^n$ mod $3$ and mod $5$. 
\end{proof}

\noindent
With this we are ready to formulate the statement that will lead us to a primality testing algorithm for $\pr_{m,n}$.
\begin{theo.}\label{ellseq}
Let $m\geq 1$ be odd and $n>1$ be such that $\pr_{m,n}:=m2^n-1$  is not divisible by $3$ or $5$ (see 
Lemma~\ref{aprcomposite}) and $4m<2^n$. Consider $t\in\Z$ such that the Jacobi symbol $\legendre{t^2+1}{\pr_{m,n}}=-1$. Take $(\tfrac{\alpha}{\beta},\tfrac{\gamma}{\delta}):=m(-1,t)\in E_t(\mathbb{Q})$ and define the sequence:
\[
x_0:=\tfrac{\alpha}{\beta},\quad x_{i+1}:=\tfrac{(x_i^2 +  t^2+1)^2}{4(x_i^3 -(t^2+1)x_i)}\quad \bmod \pr_{m,n}.
\]
Then $\pr_{m,n}$ is prime if and only if 
$x_i$ is well defined for every $0\leq i\leq n-1$ and $x_{n-1}\equiv 0 \bmod \pr_{m,n}$.
\end{theo.}
\begin{proof}
Suppose that $\pr_{m,n}$ is prime. Observe that $x_i$ equals the $x$-coordinate of the point $m2^i(-1,t)\in 
E_t(\F_{\pr_{m,n}})$. Since $2^n>4m$ we have that $n>2$, hence 
$\pr_{m,n}\equiv 3\bmod 4$. By Lemma \ref{ellcyc} using that $t^2+1\notin\F_{\pr_{m,n}}^2$ we have that $E_t(\F_{\pr_{m,n}})$ is 
cyclic. By Corollary \ref{elldiv}, the point $(-1,t)$ is not divisible by 
$2$ and $\#E_t(\F_{\pr_{m,n}})=m2^n$, hence since
$m$ is odd $m(-1,t)$ has order $2^n$ by 
Corollary \ref{ellgen}. This means that $x_{n-1}=x(2^{i-1}(\tfrac{\alpha}
{\beta},\tfrac{\gamma}{\delta}))$ equals the $x$-coordinate of the unique 
$\F_{\pr_{m,n}}$-rational point of order $2$ in $E_t(\F_{\pr_{m,n}})$, namely $(0,0)$. The fact that 
$x_i$ is well defined for $0\leq i\leq n-1$ also follows from the reasoning above.\\ \\
For the converse suppose that $\pr_{m,n}$ is not prime. Also suppose that the $x_i$ modulo $\pr_{m,n}$ 
are well defined for $0\leq i\leq n-1$
and that $x_{n-1}\equiv 0\bmod \pr_{m,n}$. Take $\ell\mid\pr_{m,n}$ the smallest prime divisor of $\pr_{m,n}$. Since $\legendre{t^2+1}{\pr_{m,n}}=-1$, 
it follows that $t^2+1\neq 0$ in $\F_\ell$.
Moreover $\ell$ is odd hence $E_t$ defines an elliptic curve over $\F_\ell$.\\
With this we have that the point $2^{n-1}m(-1,t)\in E_t(\F_\ell)$ has $x$-coordinate $0$ by assumption,
so this point equals $(0,0)$ which has order $2$.
Hence $m(-1,t)\in E_t(\F_\ell)$ has order $2^n$. 
Since $\pr_{m,n}$ is not divisible by $3$ or $5$ we have that $\ell>5$. Further $\ell\leq \sqrt{\pr_{m,n}}$ and $m(-1,t)$ generates a subgroup of order $2^{n}$ in $E_t(\F_\ell)$. Furthermore by the Hasse inequality and the fact that $\ell>5$ we have $\#E_t(\F_\ell)\leq (\sqrt{\ell} + 1)^2 < 2\ell$ , hence:
\[
2^n \leq \#E_t(\F_\ell)< 2\ell\leq 2\sqrt{\pr_{m,n}}=2\sqrt{m2^n -1}.
\]
Since $4m<2^n$ it follows that $4^n\leq 4m2^n-4<2^{2n} -4=4^n-4$ which is absurd. This contradiction shows that 
$\pr_{m,n}$ must be prime.
\end{proof}
\noindent 
The algorithm in the previous theorem uses the recursive iteration of a degree $4$ map (multiplication by $2$). In the next section we will define a primality test for other integers using a map of degree $2$ which is computationally more viable.\\ \\ 
We show an example algorithm for $\pr_{13,n}$ using Theorem \ref{ellseq}.\\
First note that $\pr_{13,n}$ is divisible by $3$ or $5$ when $n \equiv 0,1,2 \bmod 4$ by Lemma \ref{prconditions}. Hence,
the only non-trivial case to do a primality test is with the integers $13\cdot 2^{4k+3}-1$. To apply the previous theorem  we need that $4m=4\cdot 13<2^n=2^{4k+3}$, which indeed
holds for every $k\geq 1$. \\
Now we need to choose our elliptic curve $E_t$ according to Lemma \ref{ellcyc}, so we need a $t\in \Z$ such that $\legendre{t^2+1}{\pr_{13,4k+3}}=-1$. We state the following technical result as a lemma for this example:
\begin{lem.}\label{prconditions}
Let $\pr_{m,n}:=m2^n-1$ be an integer, $m\geq 1$ and $n>1$ such that one of the following conditions on $m$ and $n$ hold:
\begin{enumerate}[label=(\roman*)]
\item{$m\equiv 1\bmod 5$ and $n\equiv 3\text{ or }2\bmod 4$}
\item{$m\equiv 2\bmod 5$ and $n\equiv 2\text{ or }1\bmod 4$}
\item{$m\equiv 3\bmod 5$ and $n\equiv 0\text{ or }3\bmod 4$}
\item{$m\equiv 4\bmod 5$ and $n\equiv 1\text{ or }0\bmod 4$}
\end{enumerate}
Then $\pr_{m,n}\equiv\pm 2\bmod 5$ and therefore $5$ is not a square modulo $\pr_{m,n}$.
\end{lem.}
\begin{proof}
This is a direct calculation using the period of $2^n$ modulo $5$.
\end{proof}

The above Lemma \ref{prconditions} part (iii) allows us to use the curve $E_t$ for $t=2$ to check precisely $\pr_{13,4k+3}$ since $m=13\equiv 3\bmod 5$, hence $t^2+1=5$ is not a square modulo $\pr_{13,4k+3}$, for every $k\in \mathbb{N}$.\\With this, consider the curve $E_2$ given by $y^2 = x^3 - 5x$. \\ 
The $x$-coordinate of the point $13(-1,2)\in E_2(\mathbb{Q})$ can be computed instantly with a computer algebra software and is given by:
%
\[
x_0=-\tfrac{38867230505264472384304448711791072932034380121}{ 20648248720215880190543854206835397627372795209}.
\]
A computer program can be easily implemented to check the primality of $\pr_{13,4k+3}$ returning ``composite'' when the denominator of $x_j$ for $0\leq j \leq 4k+2$ is not a unit modulo $\pr_{13,4k+3}$ and returning ``prime'' when $x_{4k+2}\equiv 0\bmod \pr_{13,4k+3}$. A similar analysis can be done using these results for other sequences $\pr_{m,n}$.

\subsection{Primality testing using CM by $\Z[i]$ on elliptic curves}
Now we propose a primality test for integers of the form 
$\mathcal{S}_{p,n}:=p^2 16^n +1$ with $p\equiv\pm 1\bmod 
10$ prime and $p<2^n$. For 
the iteration step in the primality test we will use an 
endomorphism of an elliptic curve $E$ with $j$-invariant 
$1728$. The resulting algorithm is similar to
the one in the previous section but now using a 
degree $2$ endomorphism which is computationally better. \\

We chose these integers since $p^2 16^n + 1$ is prime 
in $\Z$ if and only if its Gaussian factor $p4^n +i$ 
(and its conjugate of course) is prime in $\Z[i]$. We 
did not choose $p^2 4^n +1$ since this integer is 
divisible by $5$ for $n$ odd. \\ \\
For the integers $\mathcal{S}_{p,n}$, is not immediate how to adapt a primality test as proposed in the previous section.\\
In the previous section we implicitly used the $\Z$-module structure of the elliptic curve $E_t$, that is, we used the action of $\Z\subset \text{End}(E_t)$ on $E_t$. Here we will use the action of $\Z[i]$ on $E$ for our primality testing purposes. \\ \\
Let $p\equiv 1\bmod 4$ and consider the elliptic curve $E/\F_p$ given by $y^2=x^3 -x$. Take $\xi\in \F_p$ such that $\xi^2=-1$. The action of $i\in\Z[i]$ on $E(\F_p)$ is defined as the ``multiplication by $i$" map:
\begin{align}\label{CM}
\begin{split}
i\colon E(\F_p) &\to E(\F_p), \\ 
(x,y)&\mapsto (-x,\xi y).
\end{split}
\end{align}
\noindent
The map $i$ is clearly an element of $\text{Aut}(E)\subset\text{End}(E)$ and $E(\F_p)$ obtains the structure of $\Z[i]$-module using the ring homomorphism
\begin{align}\label{endhomo}
\begin{split}
\Z[i] &\to \text{End}(E)\\
a+bi &\mapsto a\!+\!b\!\circ\!i.
\end{split}
\end{align}
We will use the next theorem for the rest of this section. It has an interesting story related to the last entry in Gau\ss' \textit{Tagebuch} (July $7^{th}$, 1814), discovered by Felix Klein in 1897 and published in \textit{Math. Annalen 1903} \cite{klein1903}. Gau{\ss} conjectured a way of calculating the number of points over $\F_p$ of a curve birational to the elliptic curve with $j$-invariant 1728  where $p\equiv 1\bmod 4$. Gustav Herglotz was the first to prove Gau\ss's conjecture in 1921. Here we show another elementary proof using modern language (first proved in \cite{herglotz1921} and more general in \cite{koblitzmod} and \cite{ir}). The subsequent corollary is precisely the conjecture predicted by Gau\ss.\\ For the rest of this text we will denote the composition of endomorphisms as  $ab:=a\circ b$ using the previously defined homomorphism  in (\ref{endhomo}).
\begin{theo.}\label{frobeniusdiv}
Let $p\equiv 1\bmod 4$ be a prime and let $E/\F_p$ be given by $y^2=x^3 -x$. Consider the $p^{th}$ Frobenius endomorphism $\phi_p$ and the identity map $\mathbbm{1}$. We have that $\text{End}(E)= \Z[i]$ and $\phi_p=a+bi\in \text{End}(E)$ satisfies 
$a^2+b^2=p$ and $(2+2i)\mid \left(\phi_p-\mathbbm{1}\right)$.
\end{theo.}
\begin{proof}
 We already saw that $\Z[i]\subseteq \text{End}(E)$. Since $p\not\equiv 3\bmod 4$ we have that $E$ is not supersingular (see Proposition \ref{elldiv} and take $t=0$), hence $\text{End}(E)$ is contained in the ring of integers of an imaginary quadratic field and then  $\text{End}(E)=\Z[i]$.
Moreover $p=\text{deg}(\phi_p)=a^2+b^2$.\\
 \\
Note that $2+2i\in\text{End}(E)$ is a separable map since $p\nmid \text{deg}(2+2i)=8$. We proceed to analyze its kernel since $2+2i\mid \phi_p-\mathbbm{1}$ if and only if $\text{Ker}(2+2i)\subset\text{Ker}(\phi_p-\mathbbm{1})=E(\F_p)$.\\ 
Let $P\in E(\F_p)$, we have that $(2+2i)P=(1+i)2P$, hence, if $Q\in \text{Ker}(1+i)$ is non-trivial, we have that:
\begin{equation*}
\text{Ker}(2+2i)=E[2](\F_p)\cup [2]^{-1}(Q) 
\end{equation*}
Note that $Q=(0,0)$ generates $\text{Ker}(1+i)\subset E(\F_p)$. Computing the tangent lines to $E/\F_p$ that contain $Q$, one obtains: 
\[
[2]^{-1}(Q)=\{ (\xi,\pm (1-\xi)),(-\xi,\pm (1+\xi) ) \}.
\]
Since $p\equiv 1\bmod 4$ we have that $\xi\in\F_p$ and the four points in $[2]^{-1}(Q)$ are fixed by $\phi_p$. Trivially the other four points $\{(0,0),(1,0),(-1,0),\infty\}$ in $\text{Ker}(2+2i)$ are fixed by $\phi_p$, hence $\text{Ker}(2+2i)\subset\text{Ker}(\phi_p-\mathbbm{1})$ and the result follows. 
\end{proof}
\noindent
This theorem gives a lot of information of the Frobenius endomorphism of $E$ and the precise answer to Gau{\ss}' last entry in his \textit{Tagebuch} which we will use soon. 
\begin{cor.}\label{pointscomplex}
Let $p\equiv 1\bmod 4$ and consider the elliptic curve $E/\F_p$ given by $y^2=x^3-x$,  then $\#E(\F_p)=p+1-2\alpha$ where $p=\alpha^2+\beta^2$ and if $p\equiv 1\bmod 8$ then $\alpha\equiv 1\bmod 4$, otherwise $\alpha\equiv 3\bmod 4$.
\end{cor.}
\begin{proof}
We have that $\#E(\F_p)=\deg(\phi_p-\mathbbm{1})=\deg(\alpha+\beta i -1)$, hence
\[
\#E(\F_p)=\alpha^2 +\beta^2 + 1 -2\alpha=p+1-2\alpha.
\]
Theorem \ref{frobeniusdiv} shows $2+2i\mid \alpha-1+\beta i$, hence $\alpha$ is odd and $\beta$ is even. Further we have that $8\mid \#\text{Ker}(\phi_p-\mathbbm{1})=\deg(\phi_p-\mathbbm{1})=(\alpha-1)^2+\beta^2$ since $\text{Ker}(2+2i)\subset\text{Ker}(\phi_p-\mathbbm{1})=E(\F_p)$ by the same theorem. \\
With this, since $\alpha^2 +\beta^2 =  p$ we have that 
\begin{equation}\label{tagebuch}
p+1-2\alpha \equiv 0\bmod 8.
\end{equation}
This implies the result.
\end{proof}
\noindent
We illustrate the corollary with the following example.\\
Consider the elliptic curve $E/\F_{37}$ given by $y^2=x^3-x$. We have that $37\equiv 5\bmod 8$. By the previous corollary, $37= \alpha^2 + \beta^2$, $\alpha\equiv 3\bmod 4$ and $\beta$ even, hence  $\alpha^2 +\beta^2 = 1 + 36$, $\alpha=-1$, and $\#E(\F_{37})=37+1-2(-1)=40$. \\ \\
\newcommand{\pcm}{\mathcal{S}}
The following proposition will be used to tell us the structure of $E(\F_{\pcm_{p,n}})$ as an abstract group,
given that $\pcm_{p,n}$ is prime.
\begin{prop.}\label{importanttorsion}
Let $p$ be a prime such that $p\equiv 1\bmod 8$ and $p-1$ is a square. Consider the elliptic curve $E:y^2 = x^3 - x$, then $p=(\alpha-i)(\alpha+i)$ in $\Z[i]$, $\#E(\F_p)=\alpha^2$ and $E(\F_{p})\cong \Z/(\alpha)\times\Z/(\alpha)$ as Abelian groups.
\end{prop.}
\begin{proof}
We have that $p-1=\alpha^2$ for some $\alpha\in\Z$, hence $p=(\alpha+i)(\alpha-i)$ in $\Z[i]$. 
Using Theorem \ref{frobeniusdiv}, since $p\equiv 1\bmod 8$ we have that $\#E(\F_p)=p+1-2=\alpha^2$. \\ \\
Let $\phi_p\in\text{End}(E)=\Z[i]$ be the $p^{th}$ power of Frobenius. The previous calculation shows that $\text{Tr}(\phi_p)=2$. Further $\deg \phi_p=p=\alpha^2 + 1=(\alpha+i)(\alpha-i)$, hence (after possibly changing the sign of $\alpha)$ the Frobenius endomorphism
is given by $\phi_p=\alpha i+1$. With this, if $P\in E(\F_p)$ we have that $P=\phi_p(P)=(\alpha i+1)(P)$. Hence $\alpha P=\infty$ and $P\in\ E[\alpha]\cong \Z/(\alpha)\times \Z/(\alpha)$. Since $\#E(\F_p)=\alpha^2$ we 
conclude that $E(\F_p)\cong \Z/(\alpha)\times \Z/(\alpha)$.    
\end{proof}
Now we present two corollaries that describe particular properties of the group  $E(\F_{\pcm_{p,n}})$
(again, provided $\pcm_{p,n}$ is prime). These corollaries will be used to extend the structure of $E(\F_{\pcm_{p,n}})$ to a cyclic $\Z[i]$-module in the subsequent proposition. 
\begin{cor.}\label{cardCM}
Let $\mathcal{S}_{p,n}:=p^2 16^n+1$ be prime and $n>0$. Consider the elliptic curve $E/\F_{\mathcal{S}_{p,n}}$ given by $y^2 = x^3 - x$, then $\#E(\F_{\mathcal{S}_{p,n}} )=p^2 16^n$.
\end{cor.}
\begin{proof}
Immediate from Proposition~\ref{importanttorsion}.
\end{proof}
\begin{cor.}\label{ptorrat}
Let $\pcm_{p,n}:=p^2 16^n + 1$ be prime with $p$ odd and $n>0$. Consider the elliptic curve $E:y^2=x^3-x$, then $E(\F_{\pcm_{p,n}})$ has full $p$-torsion, that is $E[p]\subset E(\F_{\pcm_{p,n}})$.
\end{cor.}
\begin{proof}
Again, this is a direct consequence of Proposition \ref{importanttorsion}.  
\end{proof}
The next proposition provides the group  $pE(\F_{\pcm_{p,n}})$ with the structure of a cyclic $\Z[i]$-module. 
\begin{prop.}\label{structCM}
Let $\pcm_{p,n}:=p^2 16^n+1$ be prime such that $p$ is odd and $n>0$. Consider the elliptic curve $E/\F_{\pcm_{p,n}}$ given by $y^2=x^3-x$, then $pE(\F_{\pcm_{p,n}})\cong \Z[i]/(1+i)^{4n})$ as cyclic $\Z[i]$-modules.
\end{prop.}
\begin{proof}
Since $n>0$ we know by Lemma \ref{cardCM} that $\#pE(\F_{\pcm_{p,n}})= 16^n$ and $pE(\F_{\pcm_{p,n}})$ is a finitely generated $\text{End}(E)$-module with $\text{End}(E)=\Z[i]$. Further, $\Z[i]$ is a PID and by the structure theorem of finitely generated modules over a PID there exists a finite sequence of ideals  $(1)\neq (z_1)\supseteq (z_2) \supseteq \ldots \supseteq (z_t)$ of $\Z[i]$, for some $t\in\mathbbm{N}$, such that
\begin{equation}\label{fgmodCM}
pE(\F_{\pcm_{p,n}})\cong \Z[i]/(z_1)\oplus\Z[i]/(z_2)\oplus\ldots\oplus \Z[i]/(z_t).
\end{equation}
This sequence of ideals implies that $z_1\mid z_2 \mid\ldots \mid z_t$. Let $\mathcal{N}\colon \Z[i]\to\Z$ be the norm map. Each direct summand has cardinality $\mathcal{N}(z_j)=z_j \bar{z_j}$ and $\mathcal{N}(z_j)\mid  16^n$.
Thus for every $j$ one concludes $\mathcal{N}(z_j)=2^{m_j}$ for for some power $m_j>0$. Hence $(z_j)=((1+i)^{m_j})\subset \Z[i]$. This implies, using $m_j>0$ for all $j$, that the $1+i$-torsion in $\bigoplus \Z[i]/(z_j)$
is isomorphic to $(\Z/2\Z)^t$.\\
Note that $\deg(1+i)=2$, hence $t=1$
and $pE(\F_{\pcm_{p,n}})\cong \Z[i]/((1+i)^{4n})$,
proving the result.       
\end{proof}
\noindent
Now we know that if $\pcm_{p,n}$ is prime, $pE(\F_{\pcm_{p,n}})$ is a cyclic $\Z[i]$-module. We need a generator of this $\Z[i]$-module to apply ideas as used in
the previous sections.\\
Similarly as in \cite{densav} we use the quadratic twist of $E$ given by the curve $E_{30}:30y^2=x^3 - x$, but now we will do a primality test on $\pcm_{p,n}$ instead of Fermat numbers. 
Assuming $\pcm_{p,n}$ is prime, the curve $E_{30}$ is isomorphic to $E:y^3=x^3-x$ over $\F_{\pcm_{p,n}}$ if and only if $\legendre{30}{\pcm_{p,n}}=1$. The following simple lemma will tell us for which $p$ the element $30$ is a square  in $\F_{\pcm_{p,m}}$.
\begin{lem.}\label{isotwist}
Let $p\equiv\pm 1\bmod 10$ and $\pcm_{p,n}:=p^2 16^n + 1$ be prime, then $E_{30}(\F_{\pcm_{p,n}})\cong E(\F_{\pcm_{p,n}})$
\end{lem.}
\begin{proof}
We have to show that $30$ is a square modulo $\pcm_{p,n}$. This is a direct application of the properties of the Legendre symbol,
using that $\pcm_{p,n}\equiv 2\bmod 5$
and $\pcm_{p,n}\equiv 2\bmod 3$ and 
$\pcm_{p,n}\equiv 1\bmod 8$; these properties imply
\begin{equation}
\legendre{30}{\pcm_{p,n}}=\legendre{5}{\pcm_{p,n}}\legendre{3}
{\pcm_{p,n}}\legendre{2}{\pcm_{p,n}}=(-1)(-1)(1)=1.
\end{equation}
\end{proof}
\noindent
The curve $E_{30}$ was chosen since the point $p(5,2)$ turns out
to be a generator of the cyclic $\Z[i]$-module $pE_{30}(\F_{\pcm_{p,n}})$. We proceed to prove this. 
\begin{lem.}\label{cmgen}
Let  $p\equiv \pm 1 \bmod 10$ be prime and consider $Q:=p(5,2)\in pE_{30}(\F_{\pcm_{p,n}})\cong \Z[i]/((1+i)^{4n})$. The point $Q$ generates the $\Z[i]$-submodule $\Z[i]/((1+i)^{4n})$ of $E_{30}(\F_{\F_{\pcm{p,n}}})$
\end{lem.}
\begin{proof}
Since $p$ is odd, we just need to show that $(5,2)$ is not in the image of $(1+i)$. This is 
the same as saying that $(1+i)(X,Y)=(X,Y)+(-X,\xi Y)= (5,2)$ has no $\F_{\pcm_{p,n}}$-rational solution.\\
We proceed to calculate $(1+i)(X,Y)$ explicitly. Consider the endomorphism $1+i\in \text{End}(E_{30})$ where $E_{30}:30y^2 = x^3 - x$. The slope between $(X,Y)$ and $i(X,Y)=(-X,\xi Y)$ is $\lambda:=\tfrac{(1-\xi)Y}{2X}$. A quick computation shows that
\begin{equation}\label{1isum}
(1+i)(X,Y)=(30\lambda^2, \lambda(X-30\lambda^2)-Y )=(\tfrac{\xi(1-{X}^2)}{2X}, - \tfrac{(1+\xi)({X}^2 +1)Y}{4{X}^2}   \big).
\end{equation}
We have that $(5,2)\in E_{30}(\F_{\pcm_{p,n}})$ is not in the image of $1+i\in\text{End}(E_{30})$ (divisible by $1+i$) if and only if the solutions of the equations below for $X$ and $Y$ are not $\F_{\pcm_{p,n}}$-rational:
\begin{equation*}
(1+i)(X,Y)=(30\lambda^2, \lambda(X-30\lambda^2)-Y )=(5,2)
\end{equation*}
If we look at the equation $30\lambda^2=5$, it means that $5$ must be a square modulo $\pcm_{p,n}$ since $30$ is by Lemma \ref{isotwist}. By the proof of that lemma we have that $\legendre{5}{\pcm_{p,n}}=-1$ so there is no such $\F_{\pcm_{p,n}}$-rational point $(X,Y)$, hence $(5,2)$ is not divisible by $1+i$. Therefore
 $p(5,2)$ generates $\Z[i]/((1+i)^{4n})$ and the generated submodule has cardinality $2^{4n}=16^n$.
\end{proof}

\noindent
Now we state the main theorem of this section. This theorem will lead us to a conclusive deterministic primality test algorithm for integers of the form $\pcm_{p,n}=p^2 16^n + 1$. 
\begin{theo.}\label{algoCMtheo}
Consider the integer $\pcm_{p,n}:=p^2 16^n +1$ where $p\equiv\pm 1\bmod 10$ is prime, and $p<2^n$. Let $E_{30}$ be the elliptic curve $30y^2=x^3-x$ and consider the point $Q:=p(5,2)\in E_{30}(\mathbb{Q}(i))$ (which is a $\Z[i]$-module).\\
Then $\pcm_{p,n}$ is prime if and only if $(1+i)^{4n-1}Q=(0,0)\bmod (p4^n+i)$. 
\end{theo.}
\begin{proof}
Suppose that the congruence holds and $\pcm_{p,n}$ is not prime. Take $k\mid \pcm_{p,n}$ the smallest prime divisor of $\pcm_{p,n}$, hence $k\leq \sqrt{p^2 16^n + 1}$. Further we have that $p^2 16^n+1\equiv 0\bmod k$ if and only if $(p4^n)^2\equiv -1\bmod k$. This means that $-1$ is a square in $\F_k$ and then $k=\pi\overline\pi$ with $\pi\in\Z[i]$ Gaussian prime.\\Now, since $\pi\mid p^2 16^n + 1=(p4^n+i)(p4^n-i)$ without loss of generality, assume that $\pi\mid p4^n +i$. Let $\mathcal{N}$ be the Gaussian norm. Since $\pi\overline\pi=k\leq \sqrt{\pcm_{p,n}}$ we have that: 
\begin{equation}\label{impineqcm}
\mathcal{N}(\pi)<\sqrt{\mathcal{N}(p4^n+i)}=\sqrt{\pcm_{p,n}}
\end{equation}
Further, the discriminant of $E_{30}$ is $(2^2\cdot 3\cdot 5)^2$ and it is easy to see that $p^2 16^n + 1$ is not divisible by $2,3$ or $5$. Hence, $E_{30}(\Z[i]/(\pi))$ defines an elliptic curve.
Furthermore, $(1+i)^{4n-1}Q\equiv (0,0)\bmod \pi$ in $E_{30}(\Z[i]/(\pi))$ if and only if $(1+i)^{4n}Q=\infty \bmod \pi$. 
This means that $Q=p(5,2)$ generates a $\Z[i]$-submodule 
of $E_{30}(\Z[i]/(\pi))$ 
of size $16^n$. With this we get the following inequalities using the Hasse inequality and the inequality in (\ref{impineqcm}): 
\begin{equation}
16^n\leq \#E_{30}(\Z[i]/(\pi))\leq (\sqrt{\mathcal{N}(\pi)}+1)^2< (\sqrt[4]{\pcm_{p,n}}+1)^2.
\end{equation}
This implies that $4^n-1< \sqrt[4]{p^2 16^n + 1}$ and then $\tfrac{(4^n-1)^4-1}{16^n}<p^2$. Since $p<2^n$ by hypothesis this implies that $\tfrac{(4^n-1)^4-1}{16^n}<4^n$, hence $0\leq n<\varepsilon$ with $\varepsilon\approx 0.91<1$ which is absurd since $n\geq 1$. \contradiction.  We conclude that $\pcm_{p,n}$ is prime. 
 \\ \\
\noindent
Suppose that $\pcm_{p,n}$ is prime, then $Q$ generates $pE_{30}(\F_{\pcm_{p,n}})\cong \Z[i]/((1+i)^{4n})$ by Lemma \ref{cmgen}. Further $(0,0)$ is the only non-trivial point in the $(1+i)$-torsion of $E_{30}(\F_{\pcm_{p,n}})$, hence $(1+i)^{4n-1}Q=(0,0)\bmod p4^n+i$ since  $E_{30}(\F_{\pcm_{p,n}})\cong E_{30}(\Z[i]/(p4^n+i))$.
\end{proof}
\noindent
The same theorem can be stated as an algorithm.
\begin{cor.}\label{algoCMfinal}
Consider the integer $\pcm_{p,n}:=p^2 16^n + 1$ such that $p$ is prime, $p\equiv\pm 1\bmod 10$ and $p<2^n$. Let $(x_0,y_0):=p(5,2)\in E_{30}(\mathbb{Q}(i))$ and consider the sequence:
\begin{equation*}
x_{j+1}= \tfrac{i(1-x_j^2)}{2x_j}\quad\bmod p4^n + i
\end{equation*}
$\pcm_{p,n}$ is prime if and only if $x_j$ is well defined for all $j<4n$ and $x_{4n-1}\equiv 0\bmod p4^n + i$
\begin{proof}
This is equivalent to Theorem \ref{algoCMtheo}. The sequence is the recursive multiplication by $(1+i)$ starting with the $x$ coordinate of the point $p(5,2)\in E_{30}$. This formula was deduced in equation (\ref{1isum}).
\end{proof}
\end{cor.}
\noindent
In order to implement the previous corollary as an algorithm to check $\pcm_{p,n}$ note that $n>\tfrac{\log(p)}{\log(2)}$ for a conclusive primality test.\\Consider the ring $\Z/(\pcm_{p,n})$. We have that $i:=p\cdot 4^n$ and $i^2 = -1$ in $\Z/(\pcm_{p,n})$. \\
Take the curve $E_{30}:30y^2 = x^3 - x$. The curve $E':y^2 = x^3 - 900x$ is isomorphic fo $E_{30}$ under the change of variables $(x,y)\mapsto (30x,900y)$.  The initial value of the iteration is $x_0$ from $(x_0,y_0)=p(5,2)$. We calculate it in $E'$ as $p(30\cdot 5,900\cdot 2)$ and its $x$ coordinate divided by $30$ will be our $x_0$. 

%
%
%
%
%
%
For the sake of completeness, we show all the primes $\pcm_{
}=p^2 16^n + 1$ with this technique such that $p\equiv\pm 1\bmod 10$, $p\leq 101$ and $\tfrac{\log p}{\log 2}<n\leq 2000$ using a GP/PARI implementation.
\begin{center}
    \begin{tabular}{| l | l | l |}
    \hline
    $p$ & $\approx \tfrac{\log p}{\log 2}$ & $n$ values where $\pcm_{p,n}=p^2 16^n + 1$ is prime $\tfrac{\log p }{\log 2}<n\leq 2000$ \\ \hline
    11 & 3.45943 & $11,21,24,57,66,80,183,197,452,1982$  \\ \hline
    19 & 4.24792 & $7,9,25,78,142,646$ \\ \hline
    29 & 4.85798 & $6,19,33,36,86,103,326,352$ \\ \hline
    31 & 4.95419 & $5,65,142,148,196,1154$  \\ \hline
    41 & 5.35755 & $12,18,48,81,113,305,620,1098$  \\ \hline
    59 & 5.88264 & $9,19,33,46,121,264,904,1365,1858$ \\ \hline
    61 & 5.93073 & $11,259,361,415,427,594$ \\ \hline
    71 & 6.14974 & $12,21,33,36,49,70,82,85,91,111,114,129,147,255$ \\ \hline
    79 & 6.30378 & $13,17,19,81,375,1027,1562,1785$ \\ \hline
    89 & 6.47573 & $39,41,47,65,71,99,299,909,1901$  \\ \hline
    101 & 6.65821& $8,202,238,1484$  \\ \hline

    \end{tabular}
\end{center}
%
%
%
%
\newcommand{\sqp}{\mathcal{T}}
\section{Primality testing using real multiplication 
on hyperelliptic Jacobians of dimension $2$}
\noindent
In this last part we propose a 
primality test 
for integers $\lambda_n:=4\cdot 5^n-1$ 
 using the Jacobian of a hyperelliptic curve of genus $2$.
This will be done similarly to what was done in the
previous section using complex multiplication on an elliptic curve.\\
This section is motivated by an open question stated by Abatzoglou, Silverberg, Sutherland and Wong  in \cite{primeframework} (Remark 4.13) asking for a primality test algorithm using higher dimensional Abelian varieties such as Jacobians of genus $2$ curves. We will use the Jacobian $\ja$ of the hyperelliptic curve $y^2 = x^5 + h$. We begin with the structure of $\text{End}(\ja)$.  \\ 
\begin{prop.}\label{endjacstruc}
Let $h\neq 0$ and $\hc:y^2 = x^5 + h$ be a hyperelliptic curve of genus $2$ over $\mathbb{Q}$. Consider the Jacobian of $\hc$ denoted by $\ja$. Then $\text{End}(\ja)=\Z[\zeta]$ where $\zeta$ is a primitive fifth root of unity.
\end{prop.}
\begin{proof}
Let $\zeta^\ast\in\text{Aut}(\hc)$ be the automorphism $\zeta^\ast(x_0,y_0)=(\zeta x_0,y_0)\in \hc$ where $\zeta$ is a primitive fifth root of unity. The action of $\zeta^\ast$ on $\hc$ is naturally extended to the Jacobian 
 which implies that $\zeta^\ast\in\text{End}(\ja)$.
 As $\zeta^\ast$ generates a subring $\cong \Z[\zeta]\subset\text{End}(\ja)$ 
 and $\ja$ is a simple Abelian variety over 
$\overline{\mathbb{Q}}$ (see \cite{PROLEGOMENA}, Chapter 15) 
, and moreover $\Z[\zeta]$ is integrally closed,
we have that 
$\text{End}(\ja)=\Z[\zeta]$.
\end{proof}
\noindent
\textbf{Remark:} We will use $\ja$ and $\text{End}(\ja)$ to test whether $\lambda_n=4\cdot 5^n-1$ is prime. Note that $3\mid\lambda_{2k}$, so we will only test $\lambda_n$ when $n$ is odd.   \\ \\
\noindent
Let $\lambda_n$ be prime. We proceed to deduce the group structure of $\ja(\F_{\lambda_n})$. First we state and prove two easy lemmas that will tell us the structure of $\ja[2](\F_{\lambda_n})$. 
\begin{lem.}\label{easyweierstrass}
Let $\hc$ be the hyperelliptic curve given by $y^2 = x^5 + h$ and let $\lambda_n:=4\cdot 5^n -1$ be prime, then there is only one $\F_{\lambda_n}$-rational point in $\hc$ of the form $(\alpha,0)$ for some $\alpha\in\F_{\lambda_n}$ 
\end{lem.}
\begin{proof}
This follows from Fermat's little theorem. We have that $5$ and $\lambda_n-1=2(2\cdot 5^n -1)$ are coprime, hence, the map $x\mapsto x^5$ is invertible over $\F_{\lambda_n}$, hence, $x^5 = -h$ has only one solution in $\F_{\lambda_n}$. 
\end{proof}
\noindent
In the situation of Lemma~\ref{easyweierstrass}, let $\alpha\in\F_{\lambda_n}$ satisfy $\alpha^5+h=0$. Observe that the zeros of $x^5 + h$ are given by $\zeta^j\alpha$ for $0 \leq j\leq 4$ and $\zeta$ a fifth root of unity. Therefore by the previous lemma $\zeta\notin\F_{\lambda_n}$. In order to deduce the structure of the $2$-torsion of $\ja$, the following lemma tells us the field extension of $\F_{\lambda_n}$ where $\zeta$ lives and this will give us directly the structure of the $2$-torsion of $\ja$. 
\noindent
\begin{lem.}\label{zetaquad}
Let $n>0$ and suppose $\lambda_n:=4\cdot 5^n -1$ is prime. We consider the field $\F_{\lambda_n}$. Let $\zeta$ be a primitive fifth root of unity, then $\zeta\in \F_{\lambda_{n}^2}$ and $\zeta^{\lambda_n}=\zeta^{-1}$.
\end{lem.}
\begin{proof}
This is immediate by observing that $\lambda_n^2\equiv 1\bmod 5$, using that the unit group of a finite field is cyclic.
\end{proof}
\noindent
With this we will deduce the structure of $\mathcal{J}[2](\F_{\lambda_n})$ in the next corollary.
\begin{cor.}\label{2torsionforprimality}
Let $n>0$ and suppose $\lambda_n$ is prime. Consider the hyperelliptic curve $\hc/\F_{\lambda_n}$ given by $y^2 = x^5 + h$. Then $\ja[2](\F_{\lambda_n})\cong \Z/(2)\times\Z/(2)$. 
\end{cor.}
\begin{proof}
We know that $\ja[2](\F_{\lambda_n})\subset\ja$ consists of divisor classes $D-2\infty$ where $D$ 
consists of  pairs of Weierstrass points of $\hc$ and $D$
is fixed under the action of the absolute Galois group of $\F_{\lambda_n}$. 
By the previous discussion we know that all the Weierstrass points of 
$\hc$ are of the form $(\zeta^j\alpha,0)$ for $0\leq j\leq 4$,
with $\alpha\in\F_{\lambda_n}$ satisfying $\alpha^5+h=0$. 
Further, only two  Weierstrass points are defined over $\F_{\lambda_n}$ 
by Lemma \ref{easyweierstrass}, namely $(\alpha,0)$ and $\infty$. 
The other four lie in a quadratic extension of $\F_{\lambda_n}$ since 
$\zeta$ lies there by Lemma \ref{zetaquad}. Let $\rho_k:=\zeta^k\alpha$ be a zero of $x^5 + h$,
then Lemma~\ref{easyweierstrass} shows that the only
conjugate of $\rho{k}$ is $\rho_{-k}$ (if $5\nmid k$).
Hence there are two pairs of conjugate Weierstrass points plus two ordered pairs of 
$\F_{\lambda_n}$-rational Weierstrass points: 
\[
\ja[2](\F_{\lambda_n})=\big\{  \{(\rho_1,0), (\rho_4,0)\}, \{(\rho_2,0), (\rho_3,0)\}, \{(\rho_0,0),\infty\}, \{\infty,\infty\}   \big\}
\]
Therefore $\ja[2](\F_{\lambda_n})\cong \Z/(2)\times\Z/(2)$. 
\end{proof}

\begin{prop.}\label{jacstruc}
Let $\hc$ be the hyperelliptic curve given by $y^2 = x^5 + h$ and suppose $\lambda_n:=4\cdot 5^n - 1$ is prime and $n>0$. Then $\#\ja(\F_{\lambda_n})=16\cdot 5^{2n}$ and $\ja(\F_{\lambda_n})\cong \Z/(\lambda_n+1)\times \Z/(\lambda_n+1)=\Z/(4\cdot 5^n)\times \Z/(4\cdot 5^n)$
\end{prop.}
\begin{proof}
First we calculate the zeta function of $\hc$. We refer to an old paper by Tate and Shafarevich \cite{cardsuper} where they proved that the numerator of the zeta function of the curve $\mathcal{C}/\F_p$ given by $y^e=x^f + \delta$ can be described explicitly when $m=\text{lcm}(e,f)|p^k+1$ for some $k$. In our case $p=\lambda_n=4\cdot 5^n -1$, $m=10$ and $k=1$. 
By \cite{cardsuper} the numerator of the zeta-function of $\hc/\F_{\lambda_n}$ is in this case given by $\lambda_n^2 T^4 + 2\lambda_n T^2 + 1$ which tells us the characteristic polynomial $\chi_\ja(T)$ of Frobenius of $\ja$ equals $ T^4 + 2\lambda_n T + \lambda_n^2=(T^2 + \lambda_n)^2$. With this $\#\ja(\F_{\lambda_n})=\chi_\ja(1)=16\cdot 5^{2n}$.\\For the structure of $\ja(\F_{\lambda_n})$, using that $\chi_\ja(T)=(T^2+\lambda_n)^2$ and $\lambda_n\equiv 3\bmod 4$, by Theorem 3.2 (iii) in \cite{onsuper}, we have that $\ja(\F_{\lambda_n})\cong \Z/(\tfrac{4\cdot 5^n}{2^a})\times \Z/(\tfrac{4\cdot 5^n}{2^b}) \times \Z/(2^{a+b})$ with $0\leq a,b\leq 2$. Further, by the previous Lemma \ref{2torsionforprimality}, $\ja[2](\F_{\lambda_n})\cong\Z/(2)\times\Z/(2)$. Hence $a=b=0$ and $\ja(\F_{\lambda_n})\cong\Z/(4\cdot 5^n)\times \Z/(4\cdot 5^n)$.  
\end{proof}

\begin{lem.}\label{lemedof}
Let $\hc$ be the hyperelliptic curve $y^2 = x^5 + h$ and take $\lambda_n:=4\cdot 5^n-1$ prime. Then $\sqrt{5}\in\text{End}_{\F_{\lambda_n}}(\ja)$
\end{lem.}
\begin{proof}
By Proposition \ref{endjacstruc} we have that $\text{End}(\ja)=\Z[\zeta]$ with $\zeta$ a primitive fifth root of unity. Using the fifth cyclotomic polynomial we have that $1+\zeta+\zeta^2 + \zeta^3+ \zeta^4=0$. Consider $\rho := \zeta +\zeta^4$, then $\rho^2 = \zeta^3 + \zeta^2 + 2=1-(\zeta+\zeta^4)=1-\rho$. With this we have that $\rho^2+\rho=1$ if and only if $4(\rho^2  + \rho)+1=5$ if and only if $(2\rho+1)^2=5$. With this $2(\zeta+\zeta^4)+1$
is a square root of $5$ in $\text{End}(\ja)$.
It is defined over $\F_{\lambda_n}$ since the Frobenius
automorphism interchanges $\zeta$ and $\zeta^4$ by Lemma~\ref{easyweierstrass}.
\end{proof}
\begin{prop.}\label{jacmodule}
Let $\lambda_n$ be prime and consider the hyperelliptic curve $\hc/\F_{\lambda_n}$ given by $y^2=x^5 + h$, then $4\ja(\F_{\lambda_n})\cong \Z[\sqrt{5}]/(\sqrt{5}^{2n})$ as $\Z[\sqrt{5}]$-modules.
\end{prop.}
\begin{proof}
By Proposition \ref{jacstruc} we have $4\ja(\F_{\lambda_n})\cong \Z/(5^n)\times\Z/(5^n)$.
This is a $\Z[\sqrt{5}]$ module with $\sqrt{5}\in \text{End}_{\F_{\lambda_n}}(\ja)$ acting
as $2(\zeta+\zeta^4)+1$. Moreover
$\sqrt{5}^{2n}$ acts trivially.
Since $\Z[\sqrt{5}]/(5^n)\cong \Z/(5^n)\times \Z/(5^n)$,
the module is necessarily cyclic since otherwise it would
contain too many elements of order $5$.
\end{proof}
\subsection{Computation of $\sqrt{5}\in\text{End}(\ja)$}
\noindent
We use the Mumford representation for elements of $\ja$ and briefly recall this here. Details and proofs of correctness and uniqueness are given in classical texts such as \cite{TLTII,CANTORJAC}. We fix our curve $\hc:y^2=x^5 +h$ and its Jacobian $\ja$.\\ 
Any point in $\ja$ is represented by a divisor $D-2\infty$
on $\hc$, with $D$ a sum of two points.
In case $D=(x_1,y_1)+(x_2,y_2)$, then define polynomials $u(x)=(x-x_1)(x-x_2)$ and $v(x)$ of degree $\leq 1$
such that $v(x_i)=y_i$. Then
\begin{equation}\label{mumf}
v(x)^2\equiv x^5 + h\bmod u(x).
\end{equation}
Note that the pair $u,v$ determines the divisor $D$.
In case $D= (x_1,y_1)+\infty$ put $u(x)=x-x_1$ and $v=y_1$,
and if $D=2\infty$ put $u=1$ and $v=0$.
So in all cases the pair $u,v$ determines $D$.\\
For the generic point $\mathfrak{g}:=(x_1,y_1)+(x_2,y_2)-2\infty\in\ja$, the coefficients of $u(x)=x^2-\alpha x+\beta$ and $v(x)=\gamma x+\delta$ are given by the symmetric functions $\alpha=x_1+x_2, \beta=x_1x_2, \gamma=\tfrac{y_1-y_2}{x_1-x_2}, \delta=\tfrac{x_2y_1-x_1y_2}{x_1-x_2}$.
The congruence (\ref{mumf}) yields defining equations for
an affine part of $\ja$ under this representation.\\ \\
With this representation, the points of $\ja$ will be denoted by $\langle u(x),v(x) \rangle$. Cantor in \cite{CANTORJAC} developed a useful algorithm to do arithmetic in $(\ja,\oplus)$ using this representation; in
fact as he explains, his method
generalizes to every hyperelliptic Jacobian of genus $g$. \\ \\
Now we show how to construct the $\sqrt{5}$ endomorphism acting on the generic point $\mathfrak{g}\in\ja(\mathbb{Q}(\sqrt{5}))\subset\ja(\mathbb{Q}(\zeta))$. Further, we will show how to deal with the exceptional case when the image of $\sqrt{5}$ corresponds to an \textit{exceptional element} (not generic) of the form $(\sigma,\rho)-\infty\in\ja$. \\
By Lemma \ref{lemedof} we have that $\eta:=\zeta+\zeta^4=\tfrac{-1+\sqrt{5}}{2}$. We know that $\zeta^i$ acts on the points of $\hc$ by multiplication on their $x$ coordinate. This action is naturally extended to $\ja$, namely 
$\zeta^i$ maps $\mathfrak{g}$ to 
$(\zeta^i x_1,y_1)+(\zeta^ i x_2,y_2)-2\infty$. With this we evaluate the 
image of the generic point under $\eta\in\text{End}(\ja)$ explicitly:
\begin{equation}\label{etadiv}
\eta(\mathfrak{g})= (\zeta x_1,y_1)+(\zeta x_2,y_2)-2\infty \oplus (\zeta^4 x_1,  y_1)+(\zeta^4 x_2,y_2)-2\infty .
\end{equation}
Let $u(x)=x^2-\alpha x+\beta$ and $v(x)=\gamma x+\delta$ be the polynomials representing the generic point $\mathfrak{g}$ of $\ja$ in Mumford representation and let $\mathfrak{G}:=\langle u(x),v(x)\rangle\in\ja$. The Mumford representation of (\ref{etadiv}) is given by the resulting divisor below which can be calculated explicitly using Cantor's addition:
\[
\mathfrak{G}_\eta:=\eta(\mathfrak{G}) = \langle x^2-\zeta\alpha x + \zeta^2\beta, \zeta^4(\gamma x +\delta) \rangle \oplus \langle x^2 - \zeta^4\alpha x + \zeta^3\beta, \zeta(\gamma x +\delta) \rangle.
\]
Then $\sqrt{5}\mathfrak{G}=2\mathfrak{G}_\eta+\mathfrak{G}$ using again Cantor's addition since $\eta=\tfrac{-1+\sqrt{5}}{2}$. The 
polynomials $u_\eta$ and $v_\eta$ defining the resulting divisor $\mathfrak{G}_\eta$ will have coefficients in $\mathbb{Z}[\eta]$ by Lemma \ref{lemedof}.\\ \\
For the case of multiplication by $\sqrt{5}$ acting on an exceptional element of the form $\hat{\mathfrak{G}}:=\langle x-x_0,y_0\rangle$, we calculate
\begin{align*}
\hat{\mathfrak{G}}_\eta:=\eta \hat{\mathfrak{G}}&=\langle x-\zeta x_0,y_0\rangle \oplus \langle x-\zeta^4 x_0,y_0\rangle\\
&=\langle x^2 -(\zeta+\zeta^4)x_0 x+x_0^2 , y_0\rangle \\
&=\langle x^2-\eta x_0 x + x_0^2,y_0	\rangle.
\end{align*}
Similarly to the previous case, we calculate the explicit formula for $\sqrt{5}$ in this exceptional case using Cantor's addition by $\sqrt{5}\hat{\mathfrak{G}}=2\hat{\mathfrak{G}}_\eta + \hat{\mathfrak{G}}$. \\ \\
The remaining case is if the resulting element of $\ja$ under $\sqrt{5}\in\text{End}(\ja)$ is exceptional (not generic), that is, $\mathfrak{D}_0\in\ja$ and $\sqrt{5}\mathfrak{D}_0=\langle x-\sigma,\rho\rangle$ or $\sqrt{5}\mathfrak{D}_0=\langle 1,0\rangle$. This can be managed in 
several ways. For example, fix a divisor $\mathfrak{D}_c\in \ja$ such that 
$\sqrt{5}\mathfrak{D}_c$ is not an \textit{exceptional element} of $\ja$.
Calculate $\mathfrak{L}:=\sqrt{5}
(\mathfrak{D}_0+\mathfrak{D}_c)$ and if $\mathfrak{L}$ results again in an 
exceptional divisor repeat this procedure with a different $\mathfrak{D}_c$. Hence 
using Cantor's addition, we obtain 
$\sqrt{5}\mathfrak{D}_0=\mathfrak{L}-\sqrt{5}\mathfrak{D}_c=\langle x-
\sigma ,\rho\rangle$. If $\mathfrak{L}=\sqrt{5}\mathfrak{D}_c$, it means that $\mathfrak{D}_0\in\text{Ker}(\sqrt{5})$ and $\sqrt{5}\mathfrak{D}_0=\langle 1,0 \rangle$ is the identity.\\
Now we are ready to formulate the main theorem-algorithm of this section.
\begin{theo.}\label{hypertest}
Let $n>1$ be an odd integer and let $\lambda_n:=4\cdot 5^n -1$.\\
Consider the hyperelliptic curve $\hc/\mathbb{Q}(\sqrt{5})$ given by $y^2\!=\!x^5+h$ with $\lambda_n\!\nmid\!h$. Suppose $\mathfrak{F}\!\in\! \ja(\mathbb{Q}(\sqrt{5}))$ is given and
and consider the sequence of divisors $\mathfrak{D}_0:=4\mathfrak{F}$, $\mathfrak{D}_{i}:=\sqrt{5}\mathfrak{D}_{i-1}=\langle u_i(x),v_i(x)\rangle$ with its coefficients reduced in $\Z[\tfrac{1+\sqrt{5}}{2}]/(2\cdot 5^{\frac{n+1}{2}}\sqrt{5}-1)\cong\Z/(\lambda_n)$.\\
If $\mathfrak{D}_j$ is well defined and $\neq \langle 1,0\rangle$ for $j\leq 2n-1$ and $\mathfrak{D}_{2n}=
\langle 1,0\rangle$  then $\lambda_n$ is prime and $\overline{\mathfrak{F}}\notin [\sqrt{5}]\ja(\F_{\lambda_n})$ for $[\sqrt{5}]\in\text{End}_{\F_{\lambda_n}}(\ja)$. 
\end{theo.}
\begin{proof}
Suppose that $\mathfrak{D}_j=\langle u_j(x),v_j(x)\rangle$ is well 
defined for $0\leq j\leq 2n-1$, $\mathfrak{D}_{2n}=\langle 1,0\rangle$ and  $\lambda_n$ is not prime. Take 
the smallest prime divisor $k\mid \lambda_n$, hence $k\leq \sqrt{4\cdot 5^n 
-1}$. Since  $2$ or $5$ 
do not divide $\lambda_n$ we have that $k\neq 2,5$. Moreover since
$k\mid \lambda_n\nmid h$ it follows that $\ja$ has good reduction at $k$. 
Finally, since in $\F_k$ we have $5=1/(4\cdot 5^{n-1})$ and $n$ is odd, it follows
that $\sqrt{5}\in\F_k$.
Consider the group $\ja(\F_k)$ which, by the argument above, is
a $\Z[\sqrt{5}]$-module. 
The assumption on $\mathfrak{D}_j$ implies that $\mathfrak{D}_j$ is well defined for $0\leq j\leq 2n-1$ in $\ja(\F_k)$. Moreover $\mathfrak{D}_0$ generates a
${\mathbb Z}[\sqrt{5}]$-submodule of $\ja(\F_k)$ of size $5^{2n}$. Further, $\#\ja(\F_k)\leq (\sqrt{k}+1)^4\leq
(\sqrt[4]{4\cdot 5^n -1}+1)^4$ by the Hasse-Weil inequality. Hence
\begin{equation*}
5^{2n}\leq \#\ja(\F_k)\leq (\sqrt[4]{4\cdot 5^n -1}+1)^4.
\end{equation*}
Since $n>1$, this inequality is false \contradiction. Therefore $\lambda_n$ is prime.\\It follows that $\overline{\mathfrak{F}}\notin [\sqrt{5}]\ja(\F_{\lambda_n})$ by the existence of the sequence of $\mathfrak{D}_n$ in the hypothesis and the cardinality of $4\ja(\F_{\lambda_n})$. 
\end{proof} 
\noindent
To implement the Theorem \ref{hypertest} as an algorithm we need a fixed $h$ and an explicit $\mathfrak{F}\in \ja$. Optionally, a proof that if $\lambda_n$ is prime then $\overline{\mathfrak{F}}\notin[\sqrt{5}]\ja(\F_{\lambda_n})$ where $[\sqrt{5}]\in\text{End}_{\F_{\lambda_n}}(\ja)$. This proof would make the above result
into an if and only if criterion. If $\lambda_n$ is prime and $\overline{\mathfrak{F}}\in [\sqrt{5}]\ja(\F_{\lambda_n})$, we will reach the identity in $4\ja(\F_{\lambda_n})$ under recursive multiplication by $\sqrt{5}$ in less than $2n$ steps.   
However we can use the previous theorem even that this $\mathfrak{F}$ is not available to find primes when the sequence in the previous theorem can be constructed. \\ \\ 
\textbf{Example:} Let $h:=10$, that is, $\HC:y^2 = x^5 + 10$. We chose the divisor $\mathfrak{F}:=(-1,3)-\infty\in\ja$ (in Mumford representation $\langle x+1,3\rangle$). 
In this case
\[
\mathfrak{D}_0=4\mathfrak{F} = \langle x^2 + \tfrac{9678206}{70644025}x + \tfrac{117106201}{70644025}, \tfrac{3088313263561}{7125156361500}x + 
    \tfrac{22033622417431}{7125156361500}\rangle
\]
regarded as polynomials in $\Z/(\lambda_n)[x]$.\\ \\
Using Theorem~\ref{hypertest} we tested primality of $\lambda_n$ for $1<n<5000$
using the above choice of $h,\mathfrak{F}$. Only for $n\in\{3,9,13,15,25,39,69,165,171,209,339,2033\}$, the integer $\lambda_n$ was
found to be prime. But there could be gaps in this list since we did not prove that our choice of $\mathfrak{F}$ is ``not divisible by'' $\sqrt{5}$ when $\lambda_n$ is prime. However, a different computation tells us that in fact this list is complete, so we conjecture that for $h=10$, our divisor $\mathfrak{F}$ turns Theorem \ref{hypertest} into a practical deterministic primality test. If our choice of $\mathfrak{F}$ turns out to be divisible by $\sqrt{5}$, to make a practical use of Theorem \ref{hypertest} and generate the sequence $\{\lambda_n\}$ of primes without gaps, future work will be to find the correct $h$ and $\mathfrak{F}$. An idea to do this is to fix $h$ and suppose that $\lambda_n$ is prime. Then one could descend through the isogeny $\sqrt{5}\in\text{End}_{\F_{\lambda_n}}(\ja)$ explicitly and find the correct $\mathfrak{F}\in \ja(\F_{\lambda_n})\setminus[\sqrt{5}]\ja(\F_{\lambda_n})$. Another idea is to solve the $\sqrt{5}$ isogeny map equated with a choice of $\mathfrak{F}\in\ja$ and show that no solutions over $\Z/(\lambda_n)$ exist for $n>1$.
\newpage
\large{\textbf{Appendix:}
\section*{$\sqrt{5}\in\text{End}_{\mathbb{Q}(\sqrt{5})}(\ja)$ for $\hc:y^2=x^5+h$ in  MAGMA}
\newsavebox\myv
\begin{lrbox}{\myv}\begin{minipage}{\textwidth}
\begin{Verbatim}[fontsize=\tiny]
// MAGMA implementation of multiplication by Square root of 5 Endomorphism for the Jacobian of the curve y^2 = x^5 + h.
// This is a Square root of 5 - rational map, so it can be modified to work with any field where 5 is a square.
// Eduardo Ruiz Duarte

Sq5 := function(D) 
F := BaseField(Parent(D));
d := Sqrt(F!5);
h := F!-Evaluate(DefiningEquation(Curve(Parent(D))), [0,0,1]);
P<A,B,C> := PolynomialRing(F,3);
R<x> := PolynomialRing(F);
J := Parent(D);

a :=F!-Coefficient(D[1],1);
b :=F! Coefficient(D[1],0);
c :=F! Coefficient(D[2],1);

vec := [a,b,c];

An := -2*A^17*B^3*h^2 - 48*A^16*B*h^3 + 27*A^15*B^4*h^2 + 2*A^14*B^3*C^2*h^2 + 644*A^14*B^2*h^3 - 204*A^13*B^5*h^2 + 48*A^13*B*C^2*h^3 - 64*A^13*h^4 + A^12*B^8*h - 25*A^12*B^4*C^2*h^2 - 3472*A^12*B^3*h^3 + 1033*A^11*B^6*h^2 - 
       632*A^11*B^2*C^2*h^3 + 32*A^11*B*h^4 - 47/2*A^10*B^9*h + 178*A^10*B^5*C^2*h^2 + 9452*A^10*B^4*h^3 - A^9*B^8*C^2*h - 3212*A^9*B^7*h^2 + 3184*A^9*B^3*C^2*h^3 + 3328*A^9*B^2*h^4 + 603/4*A^8*B^10*h - 864*A^8*B^6*C^2*h^2 - 
       13776*A^8*B^5*h^3 + 576*A^8*B*C^2*h^4 + 192*A^8*h^5 - 1/8*A^7*B^13 + 45/2*A^7*B^9*C^2*h + 5605*A^7*B^8*h^2 - 7640*A^7*B^4*C^2*h^3 - 15008*A^7*B^3*h^4 - 785/2*A^6*B^11*h + 2406*A^6*B^7*C^2*h^2 + 10928*A^6*B^6*h^3 - 
       4544*A^6*B^2*C^2*h^4 - 2432*A^6*B*h^5 - 1/2*A^5*B^14 - 253/2*A^5*B^10*C^2*h - 4940*A^5*B^9*h^2 + 8936*A^5*B^5*C^2*h^3 + 24608*A^5*B^4*h^4 + 1/8*A^4*B^13*C^2 + 444*A^4*B^12*h - 3221*A^4*B^8*C^2*h^2 - 4472*A^4*B^7*h^3 + 
       11232*A^4*B^3*C^2*h^4 + 6784*A^4*B^2*h^5 + 19/4*A^3*B^15 + 238*A^3*B^11*C^2*h + 1353*A^3*B^10*h^2 - 4984*A^3*B^6*C^2*h^3 - 13952*A^3*B^5*h^4 + 128*A^3*B*C^2*h^5 + 256*A^3*h^6 + 5/8*A^2*B^14*C^2 - 178*A^2*B^13*h + 
       1520*A^2*B^9*C^2*h^2 + 184*A^2*B^8*h^3 - 9040*A^2*B^4*C^2*h^4 - 5888*A^2*B^3*h^5 - 51/16*A*B^16 - 112*A*B^12*C^2*h + 806*A*B^11*h^2 + 1376*A*B^7*C^2*h^3 - 1520*A*B^6*h^4 - 512*A*B^2*C^2*h^5 - 512*A*B*h^6 - 9/2*B^15*C^2 - 
       129/2*B^14*h + 16*B^10*C^2*h^2 + 256*B^9*h^3 - 192*B^5*C^2*h^4 - 192*B^4*h^5;
Ad :=  A^14*B^4*h^2 - 32*A^13*B^2*h^3 - 2*A^12*B^5*h^2 + 256*A^12*h^4 + 192*A^11*B^3*h^3 - 29*A^10*B^6*h^2 - 2560*A^10*B*h^4 + 2*A^9*B^9*h + 112*A^9*B^4*h^3 + 56*A^8*B^7*h^2 + 9696*A^8*B^2*h^4 - 13/2*A^7*B^10*h - 3112*A^7*B^5*h^3 + 
       512*A^7*h^5+ 399*A^6*B^8*h^2-16736*A^6*B^3*h^4- 51/2*A^5*B^11*h + 7752*A^5*B^6*h^3 - 2560*A^5*B*h^5 + A^4*B^14 - 1438*A^4*B^9*h^2 + 11936*A^4*B^4*h^4 + 251/2*A^3*B^12*h - 6184*A^3*B^7*h^3 + 3328*A^3*B^2*h^5 - 9/2*A^2*B^15 +
       1317*A^2*B^10*h^2 - 1760*A^2*B^5*h^4 + 256*A^2*h^6 - 261/2*A*B^13*h + 304*A*B^8*h^3 - 128*A*B^3*h^5 + 81/16*B^16 - 18*B^11*h^2 + 16*B^6*h^4;

Bn := 16*A^17*B*h^3 - A^16*B^4*h^2 - 212*A^15*B^2*h^3 + 5*A^14*B^5*h^2 - 16*A^14*B*C^2*h^3 + 64*A^14*h^4 + A^13*B^4*C^2*h^2 + 1152*A^13*B^3*h^3 + 26*A^12*B^6*h^2 + 200*A^12*B^2*C^2*h^3 - 704*A^12*B*h^4 - 9/2*A^11*B^9*h - 4*A^11*B^5*C^2*h^2 - 3396*A^11*B^4*h^3 - 
      255*A^10*B^7*h^2 - 984*A^10*B^3*C^2*h^3 + 3872*A^10*B^2*h^4 + 165/4*A^9*B^10*h - 32*A^9*B^6*C^2*h^2 + 6048*A^9*B^5*h^3 - 32*A^9*B*C^2*h^4 - 192*A^9*h^5 + 9/2*A^8*B^9*C^2*h + 714*A^8*B^8*h^2 + 2452*A^8*B^4*C^2*h^3 - 13536*A^8*B^3*h^4 - 277/2*A^7*B^11*h + 
      240*A^7*B^7*C^2*h^2 - 6392*A^7*B^6*h^3 + 1024*A^7*B*h^5 - 3/8*A^6*B^14 - 147/4*A^6*B^10*C^2*h - 735*A^6*B^9*h^2 - 3368*A^6*B^5*C^2*h^3 + 28016*A^6*B^4*h^4 + 206*A^5*B^12*h - 563*A^5*B^8*C^2*h^2 + 2056*A^5*B^7*h^3 + 1472*A^5*B^3*C^2*h^4 + 640*A^5*B^2*h^5 + 
      15/8*A^4*B^15 + 193/2*A^4*B^11*C^2*h - 41*A^4*B^10*h^2 + 2668*A^4*B^6*C^2*h^3 - 28352*A^4*B^5*h^4 + 384*A^4*B*C^2*h^5 - 256*A^4*h^6 + 3/8*A^3*B^14*C^2 - 113*A^3*B^13*h + 442*A^3*B^9*C^2*h^2 + 3424*A^3*B^8*h^3 - 4624*A^3*B^4*C^2*h^4 - 5888*A^3*B^3*h^5 - 
      13/8*A^2*B^16 - 365/4*A^2*B^12*C^2*h + 257*A^2*B^11*h^2 - 960*A^2*B^7*C^2*h^3 + 9712*A^2*B^6*h^4 - 1600*A^2*B^2*C^2*h^5 - 256*A^2*B*h^6 - 3/2*A*B^15*C^2 - 23/2*A*B^14*h + 72*A*B^10*C^2*h^2 - 2488*A*B^9*h^3 + 3136*A*B^5*C^2*h^4 + 3136*A*B^4*h^5 + 9/16*B^17 + 
      21*B^13*C^2*h + 139*B^12*h^2 - 448*B^8*C^2*h^3 - 304*B^7*h^4 + 256*B^3*C^2*h^5 + 256*B^2*h^6;
Bd := A^14*B^4*h^2 - 32*A^13*B^2*h^3 - 2*A^12*B^5*h^2 + 256*A^12*h^4 + 192*A^11*B^3*h^3 - 29*A^10*B^6*h^2 - 2560*A^10*B*h^4 + 2*A^9*B^9*h + 112*A^9*B^4*h^3 + 56*A^8*B^7*h^2 + 9696*A^8*B^2*h^4 - 13/2*A^7*B^10*h - 3112*A^7*B^5*h^3 + 512*A^7*h^5 + 399*A^6*B^8*h^2 - 
      16736*A^6*B^3*h^4 - 51/2*A^5*B^11*h + 7752*A^5*B^6*h^3 - 2560*A^5*B*h^5 + A^4*B^14 - 1438*A^4*B^9*h^2 + 11936*A^4*B^4*h^4 + 251/2*A^3*B^12*h - 6184*A^3*B^7*h^3 + 3328*A^3*B^2*h^5 - 9/2*A^2*B^15 + 1317*A^2*B^10*h^2 - 1760*A^2*B^5*h^4 + 256*A^2*h^6 - 
      261/2*A*B^13*h + 304*A*B^8*h^3 - 128*A*B^3*h^5 + 81/16*B^16 - 18*B^11*h^2 + 16*B^6*h^4;

Cn := 1/5*d*A^27*B^5*C*h^3 + 308/5*d*A^26*B^3*C*h^4 - 34/5*d*A^25*B^6*C*h^3 + 1216/5*d*A^25*B*C*h^5 - 1/5*d*A^24*B^5*C^3*h^3 - 6472/5*d*A^24*B^4*C*h^4 + 132*d*A^23*B^7*C*h^3 - 308/5*d*A^23*B^3*C^3*h^4 - 24336/5*d*A^23*B^2*C*h^5 + 1/10*d*A^22*B^10*C*h^2 + 
      33/5*d*A^22*B^6*C^3*h^3 + 62412/5*d*A^22*B^5*C*h^4 - 1216/5*d*A^22*B*C^3*h^5 - 768/5*d*A^22*C*h^6 - 7223/5*d*A^21*B^8*C*h^3 + 6232/5*d*A^21*B^4*C^3*h^4 + 216864/5*d*A^21*B^3*C*h^5 + 117/20*d*A^20*B^11*C*h^2 - 626/5*d*A^20*B^7*C^3*h^3 - 73312*d*A^20*B^6*C*h^4 +
      25264/5*d*A^20*B^2*C^3*h^5 + 2240*d*A^20*B*C*h^6 - 1/10*d*A^19*B^10*C^3*h^2 + 48522/5*d*A^19*B^9*C*h^3 - 57276/5*d*A^19*B^5*C^3*h^4 - 1119488/5*d*A^19*B^4*C*h^5 + 256/5*d*A^19*C^3*h^6 - 104*d*A^18*B^12*C*h^2 + 6593/5*d*A^18*B^8*C^3*h^3 + 
      1464156/5*d*A^18*B^7*C*h^4 - 229952/5*d*A^18*B^3*C^3*h^5 - 40064/5*d*A^18*B^2*C*h^6 - 7/80*d*A^17*B^15*C*h - 119/20*d*A^17*B^11*C^3*h^2 - 218133/5*d*A^17*B^10*C*h^3 + 319524/5*d*A^17*B^6*C^3*h^4 + 3592464/5*d*A^17*B^5*C*h^5 + 3648/5*d*A^17*B*C^3*h^6 + 
      4352/5*d*A^17*C*h^7 + 7511/10*d*A^16*B^13*C*h^2 - 41887/5*d*A^16*B^9*C^3*h^3 - 4155276/5*d*A^16*B^8*C*h^4 + 1185008/5*d*A^16*B^4*C^3*h^5 - 153664/5*d*A^16*B^3*C*h^6 + 443/80*d*A^15*B^16*C*h + 1963/20*d*A^15*B^12*C^3*h^2 + 696154/5*d*A^15*B^11*C*h^3 - 
      242332*d*A^15*B^7*C^3*h^4 - 6939952/5*d*A^15*B^6*C*h^5 - 103488/5*d*A^15*B^2*C^3*h^6 - 86016/5*d*A^15*B*C*h^7 + 7/80*d*A^14*B^15*C^3*h - 16602/5*d*A^14*B^14*C*h^2 + 35178*d*A^14*B^10*C^3*h^3 + 8292816/5*d*A^14*B^9*C*h^4 - 3682624/5*d*A^14*B^5*C^3*h^5 + 
      1896768/5*d*A^14*B^4*C*h^6 - 1792/5*d*A^14*C^3*h^7 - 949/16*d*A^13*B^17*C*h - 3208/5*d*A^13*B^13*C^3*h^2 - 1603171/5*d*A^13*B^12*C*h^3 + 3239664/5*d*A^13*B^8*C^3*h^4 + 6252048/5*d*A^13*B^7*C*h^5 + 805632/5*d*A^13*B^3*C^3*h^6 + 557824/5*d*A^13*B^2*C*h^7 + 
      1/80*d*A^12*B^20*C - 109/20*d*A^12*B^16*C^3*h + 52257/5*d*A^12*B^15*C*h^2 - 518686/5*d*A^12*B^11*C^3*h^3 - 10743456/5*d*A^12*B^10*C*h^4 + 1323936*d*A^12*B^6*C^3*h^5 - 8325888/5*d*A^12*B^5*C*h^6 + 3072*d*A^12*B*C^3*h^7 + 1024/5*d*A^12*C*h^8 + 
      10271/40*d*A^11*B^18*C*h + 50723/20*d*A^11*B^14*C^3*h^2 + 2553413/5*d*A^11*B^13*C*h^3 - 5910976/5*d*A^11*B^9*C^3*h^4 + 690672*d*A^11*B^8*C*h^5 - 3209728/5*d*A^11*B^4*C^3*h^6 - 1710336/5*d*A^11*B^3*C*h^7 + 99/320*d*A^10*B^21*C + 2151/40*d*A^10*B^17*C^3*h - 
      93387/4*d*A^10*B^16*C*h^2 + 1078094/5*d*A^10*B^12*C^3*h^3 + 6953208/5*d*A^10*B^11*C*h^4 - 5195616/5*d*A^10*B^7*C^3*h^5 + 21556736/5*d*A^10*B^6*C*h^6 - 768/5*d*A^10*B^2*C^3*h^7 + 44032/5*d*A^10*B*C*h^8 - 1/80*d*A^9*B^20*C^3 - 20489/40*d*A^9*B^19*C*h - 
      143063/20*d*A^9*B^15*C^3*h^2 - 2503679/5*d*A^9*B^14*C*h^3 + 6470872/5*d*A^9*B^10*C^3*h^4 - 16365984/5*d*A^9*B^9*C*h^5 + 7713856/5*d*A^9*B^5*C^3*h^6 + 3501312/5*d*A^9*B^4*C*h^7 + 2048/5*d*A^9*C^3*h^8 - 121/32*d*A^8*B^22*C - 15857/80*d*A^8*B^18*C^3*h + 
      677709/20*d*A^8*B^17*C*h^2 - 1454286/5*d*A^8*B^13*C^3*h^3 + 1270708/5*d*A^8*B^12*C*h^4 - 2742496/5*d*A^8*B^8*C^3*h^5 - 34085632/5*d*A^8*B^7*C*h^6 - 243712/5*d*A^8*B^3*C^3*h^7 - 482304/5*d*A^8*B^2*C*h^8 - 103/320*d*A^7*B^21*C^3 + 31121/80*d*A^7*B^20*C*h + 
      69012/5*d*A^7*B^16*C^3*h^2 + 1123937/5*d*A^7*B^15*C*h^3 - 2627588/5*d*A^7*B^11*C^3*h^4 + 19237056/5*d*A^7*B^10*C*h^5 - 2344128*d*A^7*B^6*C^3*h^6 - 1323776*d*A^7*B^5*C*h^7 - 37888/5*d*A^7*B*C^3*h^8 - 4096*d*A^7*C*h^9 + 479/40*d*A^6*B^23*C + 
      22833/80*d*A^6*B^19*C^3*h - 119693/4*d*A^6*B^18*C*h^2 + 1015729/5*d*A^6*B^14*C^3*h^3 - 5018208/5*d*A^6*B^13*C*h^4 + 8793264/5*d*A^6*B^9*C^3*h^5 + 30290176/5*d*A^6*B^8*C*h^6 + 437248/5*d*A^6*B^4*C^3*h^7 + 1115136/5*d*A^6*B^3*C*h^8 + 1111/320*d*A^5*B^22*C^3 + 
      23283/80*d*A^5*B^21*C*h - 74901/5*d*A^5*B^17*C^3*h^2 + 19771*d*A^5*B^16*C*h^3 - 358456*d*A^5*B^12*C^3*h^4 - 10588992/5*d*A^5*B^11*C*h^5 + 10547904/5*d*A^5*B^7*C^3*h^6 + 9256448/5*d*A^5*B^6*C*h^7 + 169984/5*d*A^5*B^2*C^3*h^8 + 172032/5*d*A^5*B*C*h^9 - 
      471/40*d*A^4*B^24*C - 6883/80*d*A^4*B^20*C^3*h + 166009/10*d*A^4*B^19*C*h^2 - 109313/5*d*A^4*B^15*C^3*h^3 + 2098348/5*d*A^4*B^14*C*h^4 - 5809616/5*d*A^4*B^10*C^3*h^5 - 11683328/5*d*A^4*B^9*C*h^6 + 61696*d*A^4*B^5*C^3*h^7 - 233472/5*d*A^4*B^4*C*h^8 + 
      4096/5*d*A^4*C^3*h^9 - 1351/160*d*A^3*B^23*C^3 - 76881/80*d*A^3*B^22*C*h + 35323/5*d*A^3*B^18*C^3*h^2 - 151804/5*d*A^3*B^17*C*h^3 + 1657112/5*d*A^3*B^13*C^3*h^4 + 2067584/5*d*A^3*B^12*C*h^5 - 4265792/5*d*A^3*B^8*C^3*h^6 - 5609984/5*d*A^3*B^7*C*h^7 - 
      191488/5*d*A^3*B^3*C^3*h^8 - 249856/5*d*A^3*B^2*C*h^9 + 1881/320*d*A^2*B^25*C - 5421/40*d*A^2*B^21*C^3*h - 21057/4*d*A^2*B^20*C*h^2 - 29962*d*A^2*B^16*C^3*h^3 + 318012/5*d*A^2*B^15*C*h^4 + 839232/5*d*A^2*B^11*C^3*h^5 - 1024/5*d*A^2*B^10*C*h^6 - 
      821504/5*d*A^2*B^6*C^3*h^7 - 958464/5*d*A^2*B^5*C*h^8 - 16384/5*d*A^2*B*C^3*h^9 - 16384/5*d*A^2*C*h^10 + 81/40*d*A*B^24*C^3 + 48717/80*d*A*B^23*C*h - 8412/5*d*A*B^19*C^3*h^2 - 72417/5*d*A*B^18*C*h^3 + 38432/5*d*A*B^14*C^3*h^4 + 214832/5*d*A*B^13*C*h^5 - 
      73472/5*d*A*B^9*C^3*h^6 - 115712/5*d*A*B^8*C*h^7 - 53248/5*d*A*B^4*C^3*h^8 - 53248/5*d*A*B^3*C*h^9 - 513/320*d*B^26*C + 2187/10*d*B^22*C^3*h + 2061/4*d*B^21*C*h^2 - 8552/5*d*B^17*C^3*h^3 - 12948/5*d*B^16*C*h^4 + 22528/5*d*B^12*C^3*h^5 + 19776/5*d*B^11*C*h^6 - 
      11264/5*d*B^7*C^3*h^7 - 11264/5*d*B^6*C*h^8;
Cd := A^25*B^6*h^3 - 48*A^24*B^4*h^4 - 6*A^23*B^7*h^3 + 768*A^23*B^2*h^5 + 480*A^22*B^5*h^4 - 4096*A^22*h^6 - 32*A^21*B^8*h^3 - 10752*A^21*B^3*h^5 + 3*A^20*B^11*h^2 - 456*A^20*B^6*h^4 + 73728*A^20*B*h^6 + 203*A^19*B^9*h^3 + 51408*A^19*B^4*h^5 - 87/4*A^18*B^12*h^2 - 
      11244*A^18*B^7*h^4 - 573952*A^18*B^2*h^6 + 966*A^17*B^10*h^3 - 26256*A^17*B^5*h^5 - 12288*A^17*h^7 - 129/4*A^16*B^13*h^2 + 42660*A^16*B^8*h^4 + 2522624*A^16*B^3*h^6 + 3*A^15*B^16*h - 7822*A^15*B^11*h^3 - 723360*A^15*B^6*h^5 + 159744*A^15*B*h^7 + 
      1239/2*A^14*B^14*h^2 + 52464*A^14*B^9*h^4 - 6851328*A^14*B^4*h^6 - 51/2*A^13*B^17*h + 5388*A^13*B^12*h^3 + 3342672*A^13*B^7*h^5 - 847104*A^13*B^2*h^7 - 3783/4*A^12*B^15*h^2 - 704364*A^12*B^10*h^4 + 11822720*A^12*B^5*h^6 - 12288*A^12*h^8 + 579/16*A^11*B^18*h + 
      82158*A^11*B^13*h^3 - 7249392*A^11*B^8*h^5 + 2325504*A^11*B^3*h^7 + A^10*B^21 - 12099/2*A^10*B^16*h^2 + 1972080*A^10*B^11*h^4 - 12835328*A^10*B^6*h^6 + 98304*A^10*B*h^8 + 1227/4*A^9*B^19*h - 300951*A^9*B^14*h^3 + 8637264*A^9*B^9*h^5 - 3439872*A^9*B^4*h^7 - 
      39/4*A^8*B^22 + 110025/4*A^8*B^17*h^2 - 2646732*A^8*B^12*h^4 + 8383872*A^8*B^7*h^6 - 276480*A^8*B^2*h^8 - 23337/16*A^7*B^20*h + 454166*A^7*B^15*h^3 - 5436528*A^7*B^10*h^5 + 2550528*A^7*B^5*h^7 - 4096*A^7*h^9 + 583/16*A^6*B^23 - 45762*A^6*B^18*h^2 + 
      1722120*A^6*B^13*h^4 - 3004352*A^6*B^8*h^6 + 304128*A^6*B^3*h^8 + 41163/16*A^5*B^21*h - 310965*A^5*B^16*h^3 + 1514304*A^5*B^11*h^5 - 748800*A^5*B^6*h^7 + 12288*A^5*B*h^9 - 4077/64*A^4*B^24 + 33018*A^4*B^19*h^2 - 424848*A^4*B^14*h^4 + 483328*A^4*B^9*h^6 - 
      89088*A^4*B^4*h^8 - 15525/8*A^3*B^22*h + 72254*A^3*B^17*h^3 - 143088*A^3*B^12*h^5 + 64512*A^3*B^7*h^7 - 4096*A^3*B^2*h^9 + 3159/64*A^2*B^25 - 14985/2*A^2*B^20*h^2 + 22344*A^2*B^15*h^4 - 18048*A^2*B^10*h^6 + 3072*A^2*B^5*h^8 + 7047/16*A*B^23*h - 1809*A*B^18*h^3
      + 2256*A*B^13*h^5 - 768*A*B^8*h^7 - 729/64*B^26 + 243/4*B^21*h^2 - 108*B^16*h^4 + 64*B^11*h^6;

an := Evaluate(An,vec);
ad := Evaluate(Ad,vec);
bn := Evaluate(Bn,vec);
bd := Evaluate(Bd,vec);
cn := Evaluate(Cn,vec);
cd := Evaluate(Cd,vec);

aa := an/ad;
bb := bn/bd;
cc := cn/cd;

// We can dedice D by A,B,C in <x^2 + Ax + B, Cx + D>
sigma1 :=(1/2)*(2*h - aa^5 - 5*aa*bb^2 + 5*aa^3*bb - cc^2 *(aa^2 - 4*bb));
DD := cc*(h  - aa^5 - 5*aa*bb^2 + 5*aa^3*bb + sigma1  +  aa*(aa^2-3*bb)*(aa^2-bb))/(aa^4 -bb*(3*aa^2 - bb));
Rdiv := J![x^2 + aa*x + bb, cc*x+DD ];
return Rdiv;
end function;

// We test with H:y^2 = x^5 + 10 with the ordinary divisor,  D:=2*<x+1,3>
// We calculate the multiplication by d:=SQ5 twice and compare it with 5*D

Q<d> := QuadraticField(5);
P<X> := PolynomialRing(Q);
J := Jacobian(HyperellipticCurve(X^5 + 10));
D0 := J![X+1,3];
printf "Divisor D to test\n";
D := 2*D0;
D;
printf "Square root of 5 acting on D\n";
Ds5 := Sq5(D);
Ds5; 
printf "Square root of 5 acting two times on D\n";
D5 := Sq5(Ds5);
D5;
printf "Is 5*D the same as twice the action of square root of 5 on D? %o\n", D5 eq 5*D;

\end{Verbatim}
\end{minipage}\end{lrbox}
\resizebox{0.40\textwidth}{!}{\usebox\myv}
\section*{Acknowledgements}
\noindent
This work was supported by the National Council of Science and Technology (CONACyT M\'{e}xico) through the agreement CVU-440153.\\I would like to thank my PhD supervisor Prof. Jaap Top and colleagues Ane S.I. Anema, Max Kronberg \& Marc P. Noordman for their useful comments in earlier versions of this document. Also thanks to Prof. Michael Stoll for his useful advice on MAGMA to get explicit formulas for the genus $2$ case.    
\newpage




\bibliographystyle{model1-num-names}







\end{document}